\documentclass{article}
\usepackage{maa-monthly}
\usepackage[colorlinks=true]{hyperref}

\theoremstyle{theorem}
\newtheorem{theorem}{Theorem}
\newtheorem{algorithm}[theorem]{Algorithm}
\newtheorem{lemma}[theorem]{Lemma}
\newtheorem{corollary}[theorem]{Corollary}
 
\theoremstyle{definition}
\newtheorem{definition}[theorem]{Definition}
\newtheorem{remark}[theorem]{Remark}

\begin{document}

\title{Zaremba's Conjecture for Geometric Sequences: An Algorithm}
\markright{Zaremba's Conjecture for Geometric Sequences}
\author{Elias Dubno}

\maketitle

\begin{abstract}
Even though Zaremba's conjecture remains open, Bourgain and Kontorovich solved the problem for a full density subset. Nevertheless, there are only a handful of \mbox{explicit sequences} known to satisfy the strong version of the conjecture, all of which were obtained using essentially the same algorithm. In this note, we provide a refined algorithm using the folding lemma for continued fractions, which both generalizes and improves on the old one. As a result, we uncover new examples that fulfill the strong version of Zaremba’s conjecture.
\end{abstract}

\noindent
\section{Introduction.} In 1971, S. K. Zaremba \cite{Zaremba} conjectured that every positive integer can be obtained as the denominator of a continued fraction with all partial quotients bounded by some absolute constant. More explicitly, he postulated the existence of some uniform bound $A\in\mathbb{N}$ such that for any integer $d\geq2$, there is a reduced fraction $\frac{b}{d}=[0,a_1,a_2,\dots,a_n]$ satisfying $\max{a_j}\leq A$.

\begin{definition}
    We say that such an integer $d$ is \emph{$A$-Zaremba}.
\end{definition}

Here, $\frac{b}{d}=[0,a_1,a_2\dots,a_n]$ denotes the \emph{simple continued fraction expansion} $$ \frac{b}{d}=\cfrac{1}{a_1+\cfrac{1}{a_2+\cfrac{1}{\ddots+\cfrac{1}{a_n}}}},$$ where the integers $a_j\geq1$ are called the \emph{partial quotients} of $\frac{b}{d}$.  Every rational number has a finite continued fraction expansion which is unique if we require $a_n>1$.

In its strongest form, Zaremba's conjecture even specifies the value of $A$, namely $A=5$. The simple explanation is that for $A=4$ there are known counterexamples ($d=6,54,150$)%\cite{oeis}
, while there are none for $A=5$. Hence we may reformulate Zaremba's strong conjecture as follows: every integer $d\geq2$ is $5$-Zaremba.

To this day, the conjecture remains unproven. Nevertheless, a significant breakthrough emerged in 2014, when Bourgain and Kontorovich \cite{bourgainkontorovich} confirmed the conjecture with $A=50$ through an intricate usage of the circle method for a subset of \emph{full density}. 

\begin{theorem}[\cite{bourgainkontorovich}, Theorem 1.2]
    Let $A=50$. As $N\longrightarrow \infty$, we have $$\frac{\# \{d\leq N \mid d \text{ is $A$-Zaremba} \}}{N} \longrightarrow 1.$$
\end{theorem}
By refining their methods, Huang \cite{huang} strengthened their result, showing that one may also take $A=5$.

\section{Explicit Constructions.} In addition to its inherent interest, the problem of bounding partial quotients plays an important role in the theory of good lattice points for numerical integration or generating pseudo-random numbers, which is also the context in which Zaremba originally studied the matter. (We refer to the survey \cite{kontorovich} for further elaboration.) Due to those reasons, there exists an interest in explicit examples.

The first concrete examples were exhibited by Niederreiter \cite{Niederreiter1986} who showed that powers of 2, 3 and 5 satisfy Zaremba's strong conjecture. A similar result for powers of 6 was obtained in 2002 by Yodphotong and Laohakosol \cite{YodphotongLaohakosol2002}. 

The proofs of these results are nearly identical and can be viewed as specific instances of the same algorithm. This algorithm can be roughly summarized as follows:

\begin{algorithm}\label{ALG:OLD}
    If $d$ satisfies some mild conditions, then all powers of $d$ are $A$-Zaremba with $A=d-1$.
\end{algorithm}

We will go into more details on this algorithm in the next section. 

As an application of this algorithm, Niederreiter (and later Yodphotong and Laohakosol) found concrete examples of geometric sequences satisfying Zaremba's strong conjecture.

\begin{corollary}
    Powers of $2$, $3$, $5$ and $6$ are $5$-Zaremba.
\end{corollary}

Our main result in this note is a new and refined algorithm which generalizes and improves on the old bound.

\begin{algorithm}\label{ALG:NEW}
    Let $d=x^2y$, where $x,y$ are positive integers with $xy\geq4$. If $d$ as well as $xd$ are $(xy-1)$-Zaremba such that the first and last partial quotients satisfy 
\begin{equation}\label{eq:conditionsstar}
2\leq a_1,a_n\leq xy-2,
\end{equation}
then all powers of $d$ are $A$-Zaremba with $A=xy-1$.
\end{algorithm}

\begin{remark}
   For each power of $d$ as denominator, the algorithms imply the existence of a suitable numerator such that the fraction is $A$-Zaremba. However, they do not provide any information on \emph{how many} such numerators exist. Quantitative results of this kind have for example been obtained by Kan and Krotkova in \cite[Theorems 4 \& 5]{KanKrotkova}.
\end{remark}

Before we further investigate our refined algorithm, we present new examples satisfying Zaremba's strong conjecture as an immediate corollary.

\begin{corollary}
    Powers of $12$ and $18$ are $5$-Zaremba. 
\end{corollary}
\begin{proof}
First, consider $x=2$ and $y=3$. Then $d=x^2y=12$ and $A=xy-1=5$. It suffices to check that
\begin{align*}
    \frac{5}{12} = [0,2,2,2] ~ \text{ and } ~ \frac{5}{24} = [0,4,1,4]  
\end{align*}
satisfy \eqref{eq:conditionsstar} to deduce that all powers of 12 are 5-Zaremba.

Similarly, taking $x=3$, $y=2$ and observing
\begin{align*}
    \frac{5}{18} = [0,3,1,1,2] ~ \text{ and } ~ \frac{17}{54} = [0,3,5,1,2]
\end{align*}
is enough to conclude that all powers of 18 are 5-Zaremba.
\end{proof}

Note that the value of $A$ obtained with our new algorithm is given by $xy-1$, compared to $x^2y-1$ with the old one. Therefore, for numbers $d=x^2y$ that are non-squarefree, i.e., with $x\neq1$, our new bound improves on the old one, whereas for $x=1$ we recover the bound from the old algorithm. 

\section{Idea of Algorithm \ref{ALG:OLD}.}

At the heart of both the old and the new algorithm lies a famous construction commonly known as the \emph{folding lemma} \cite{mendesfrance, vanderPoortenShallit}, which is of significance on its own.

\begin{theorem}[Folding Lemma]\label{thm:cfmatrices}
    Let $\frac{b}{d}=[0,a_1,a_2,\dots,a_n]$ be a reduced fraction, and let $z\geq1$ be a positive integer. Then the fraction $$\frac{b}{d}+\frac{(-1)^n}{zd^2}=\frac{zbd+(-1)^n}{zd^2}$$ is reduced and its continued fraction expansion is given by
\begin{align*}
    \frac{zbd+(-1)^n}{zd^2} = [0,& a_1,a_2,\dots,a_{n-1},a_n,z-1,1,a_n-1,a_{n-1},a_{n-2},\dots,a_1], 
\end{align*}
where we use the convention that $[\dots,a,0,a',\dots]=[\dots,a+a',\dots]$.
\end{theorem}

The following classical matrix representation helps prove the folding lemma. 

\begin{theorem}\label{thm}
    For $x=[0,a_1,a_2,\dots]$ and for all $m\geq1$, we have
    $$\begin{pmatrix}
        0 & 1 \\ 1 & 0
    \end{pmatrix}\begin{pmatrix}
        a_1 & 1 \\ 1 & 0
    \end{pmatrix}\begin{pmatrix}
        a_2 & 1 \\ 1 & 0
    \end{pmatrix}\dots\begin{pmatrix}
        a_m & 1 \\ 1 & 0
    \end{pmatrix}=\begin{pmatrix}
        p_m & p_{m-1} \\ q_m & q_{m-1}
    \end{pmatrix},$$
    where $\frac{p_m}{q_m}=[0,a_1,\dots,a_m]$ denotes the $m$-th convergent to $x$.
\end{theorem}

Since taking transposes in the matrix product above results in flipping the order of the partial quotients, the folding lemma now follows directly from
\begin{align*}
    \begin{pmatrix}
        p_n & p_{n-1} \\ q_n & q_{n-1}
    \end{pmatrix}\begin{pmatrix}
        z-1 & 1 \\ 1 & 0
    \end{pmatrix}\begin{pmatrix}
        1 & 1 \\ 1 & 0
    \end{pmatrix}& \begin{pmatrix}
        q_n-q_{n-1} & p_n-p_{n-1} \\ q_{n-1} & p_{n-1}
    \end{pmatrix} \\
    & =  \begin{pmatrix}
        zp_nq_n+(-1^n) & * \\ zq_n^2 & * \end{pmatrix}.
\end{align*}

\begin{definition}
     We say that we obtain the reduced fraction $\frac{zbd+(-1)^n}{zd^2}$ by applying a \emph{$z$-fold} to $\frac{b}{d}$, and we write $$\frac{b}{d} \overset{z} \longrightarrow \frac{*}{zd^2}.$$
\end{definition}

We omit writing down the explicit numerator when applying the folding lemma, as we are only interested in the denominator and we do not have to keep track of the numerator when we want to apply it successively: the new denominator does not depend on the previous numerator. This is what makes the folding lemma so handy! 

There is one technicality we have to be aware of: as mentioned in the introduction, any rational number has a unique continued fraction expansion if we require the last partial quotient to satisfy $a_n>1$. To avoid any problems and ensure that this condition is preserved, we will apply the folding lemma only to continued fractions with $a_1>1$.

Using this convention, the multiset of partial quotients of $\frac{*}{zd^2}$ after one $z$-fold is given by $\{a_1,\dots,a_n,z-1,1,a_n-1\}$ if $z>1$ and $\{a_1,\dots,a_{n-1},a_n+1,a_n-1\}$ if $z=1$. So imposing the right conditions on $z$ ensures that if $d$ is $A$-Zaremba then so is $zd^2$.

But how can one actually utilize the folding lemma to produce sequences of $A$-Zaremba integers? For instance, suppose that $\frac{b}{d}=[0,a_1,\dots,a_n]$ is $A$-Zaremba and assume that $1<a_1,a_n<A$. Then we can apply repeated 1-folds to $\frac{b}{d}$: $$\frac{b}{d} \overset{1} \longrightarrow \frac{*}{d^2} \overset{1} \longrightarrow \frac{*}{d^{2^2}} \overset{1} \longrightarrow \frac{*}{d^{2^3}} \overset{1} \longrightarrow \dots.$$ This way we find that all elements of the sequence $\bigl(d^{2^k}\bigr)_{k\geq0}$ are $A$-Zaremba. (See also \cite{komatsu}.) Some may already be happy with that, but others may note that such sequences are of \emph{doubly}-exponential growth, and wonder if one could do better. And indeed, \mbox{Niederreiter} showed that we can! Let us sketch the idea of Algorithm \ref{ALG:OLD}.

\begin{proof}[Sketch]
Let $A=d-1$ and suppose that $\frac{*}{d^k}$ is $A$-Zaremba such that the first and last partial quotients satisfy $2\leq a_1,a_n<A$. Observe that we can apply arbitrarily many 1-folds and $d$-folds to $\frac{*}{d^k}$ and always obtain a new power of $d$ as denominator which is also $A$-Zaremba.

Now all we have to do is check the first few cases (i.e., $d^k$ for small values of $k$) by hand, and then finish off using a simple inductive argument as follows: take any integer $k\geq2$. Assume that $d^m$ is $A$-Zaremba for all $m<k$ (and satisfies the mild conditions on the first and last partial quotient). If $k=2m$ is even, then we apply a 1-fold to $\frac{*}{d^m}$ to show that $d^k$ is $A$-Zaremba. Similarly if $k=2m+1$ is odd, we apply a $d$-fold to $\frac{*}{d^m}$. 
\end{proof}

So if we want to find sequences of exponential growth that satisfy Zaremba's strong conjecture, the strategy from Algorithm \ref{ALG:OLD} reaches its limits at powers of $6$. 

\section{Proof of Algorithm \ref{ALG:NEW}.} Consider a positive integer $d$ of the form $d=x^2y$. Then we can apply an $x$-fold followed by a $y$-fold to any $\frac{*}{d^k}$ and obtain a new power of $d$ in the denominator: $$\frac{*}{d^k} \overset{x} \longrightarrow \frac{*}{xd^{2k}} \overset{y} \longrightarrow \frac{*}{y(xd^{2k})^2}=\frac{*}{yx^2 d^{4k}}=\frac{*}{d^{4k+1}}.$$ 

We use this observation in combination with a more complicated inductive argument to demonstrate Algorithm \ref{ALG:NEW}. 

\begin{proof}[Proof of Algorithm \ref{ALG:NEW}]
    By assumption, $\frac{*}{d}$ satisfies \eqref{eq:conditionsstar}. Reading off the folding lemma we note that for any $z\in\{1,\dots,xy\}$, the resulting fraction after performing a $z$-fold $$\frac{*}{d} \overset{z} \longrightarrow \frac{*}{zd^2}$$
    also satisfies \eqref{eq:conditionsstar}, and hence $zd^2$ is $(xy-1)$-Zaremba as well. 
    
    We shall prove by strong induction that $d^k$ satisfies \eqref{eq:conditionsstar} for all $k\geq1$. First, we show that this holds for a special class of $k$'s.

\begin{lemma}\label{lem:basecase}
    Let $d=x^2y$ be as in Algorithm \ref{ALG:NEW}. All powers $d^k$, where $k$ is of the form $k=2^j-1$, satisfy \eqref{eq:conditionsstar} for all $j\geq1$.
\end{lemma}

\begin{proof}
    For $j=1=k$, this is true by assumption, so suppose that $j\geq2$. 

    Since $\frac{*}{xd}$ satisfies \eqref{eq:conditionsstar} by assumption as well, we can apply $(j-2)$ $xy$-folds followed by a single $y$-fold to $\frac{*}{xd}$. As explained above, the resulting fraction will satisfy \eqref{eq:conditionsstar}, and as we show now, its denominator is given by  $d^{2^j-1}$. 

    \begin{align*}
        \frac{*}{xd} \overset{xy} \longrightarrow \frac{*}{(xy)(xd)^2} & \overset{xy} \longrightarrow \frac{*}{(xy)^3(xd)^4} \overset{xy} \longrightarrow \dots \overset{xy} \longrightarrow \frac{*}{(xy)^{2^{j-2}-1}(xd)^{2^{j-2}}} \\
        & \overset{y} \longrightarrow \frac{*}{y(xy)^{2^{j-1}-2}(xd)^{2^{j-1}}}.
    \end{align*}

    We can rewrite the last denominator as $$ x^{2^{j-1}-2+2^{j-1}}y^{1+2^{j-1}-2}d^{2^{j-1}} = (x^2y)^{2^{j-1}-1}d^{2^{j-1}} = d^{2^{j-1}-1}d^{2^{j-1}} = d^{2^j-1}.$$

    Thus all powers of the form $d^k$, where $k=2^j-1$, satisfy \eqref{eq:conditionsstar} for all $j\geq1$. \end{proof}

Lemma \ref{lem:basecase} takes care of the situation $k=2^j-1$ for $j\geq1$. We now want to perform the induction step. For this, take any $k\in\mathbf{N}$. Assume that $k\neq2^j-1$, $j\geq1$, and assume that $k$ is odd. (We will deal with even $k$ later.) Then either $k \equiv 1 \mod 4$ or $k \equiv 3 \mod 4$. In the latter case, we either have $k \equiv 3 \mod 8$ or $k \equiv 7 \mod 8$. Continuing in this manner, we eventually find some $j\geq1$ such that $$k \equiv 2^j-1 \mod 2^{j+1}.$$   

Then $m := \frac{k-(2^j-1)}{2^{j+1}}$ is a non-negative integer, and by Lemma \ref{lem:basecase} we may assume $m\geq1$. 

Suppose that the claim holds for $d^m$, i.e., $\frac{*}{d^m}$ satisfies \eqref{eq:conditionsstar}. We can then apply one $x$-fold followed by $(j-1)$ $xy$-folds and one $y$-fold to $\frac{*}{d^m}$, resulting in a new reduced fraction satisfying \eqref{eq:conditionsstar} and with denominator $d^k$. The calculations are similar to the ones performed in the proof of Lemma \ref{lem:basecase}.

\begin{align*}
    \frac{*}{d^m}  \overset{x} \longrightarrow \frac{*}{x(d^m)^2} \overset{xy} \longrightarrow & \frac{*}{(xy)x^2(d^m)^4} \overset{xy} \longrightarrow \dots \overset{xy} \longrightarrow \frac{*}{(xy)^{2^{j-1}-1}x^{2^{j-1}}(d^m)^{2^j}} \\
    & \overset{y} \longrightarrow \frac{*}{y(xy)^{2^{j}-2}x^{2^{j}}(d^m)^{2^{j+1}}}.
\end{align*}

The last denominator can be rewritten as 
$$ x^{2^j-2+2^j}y^{1+2^j-2}d^{2^{j+1}m} = (x^2y)^{2^j-1}d^{2^{j+1}m} = d^{2^j-1+2^{j+1}m} = d^k. $$
Finally, consider $d^k$ for some even power $k$. Then $\frac{k}{2}=m\in\mathbf{N}$, and if $\frac{*}{d^m}$ satisfies \eqref{eq:conditionsstar}, then we can apply a 1-fold to get
$\frac{*}{d^m} \overset{1} \longrightarrow \frac{*}{(d^m)^2}=\frac{*}{d^{2m}}=\frac{*}{d^k}$. In both cases, it follows that $d^k$ satisfies \eqref{eq:conditionsstar}, concluding the proof of Algorithm \ref{ALG:NEW}. 
\end{proof}

\begin{acknowledgment}{Acknowledgment.}
I would like to thank Anna Theorin Johansson as well as my advisor, Claire Burrin, for many engaging discussions and helpful feedback. Moreover, I am very grateful to the reviewers for their constructive comments. This work was supported by Swiss National Science Foundation grant PR00P2\_201557.
\end{acknowledgment}

\begin{biog}
\item
\begin{affil}
Institute of Mathematics, University of Zurich, Winterthurerstrasse 190, 8057 Zurich, Switzerland\\
elias.dubno@math.uzh.ch
\end{affil}
\end{biog}
\vfill\eject


\begin{thebibliography}{10}

\bibitem{bourgainkontorovich} Bourgain, J., Kontorovich, A. (2014). On Zaremba's conjecture. \textit{Ann. of Math. (2)} 180(1): 137--196. doi.org/10.4007/annals.2014.180.1.3

\bibitem{huang} Huang, S. (2015). An Improvement to Zaremba’s Conjecture. \textit{Geom. Funct. Anal.} 25(3): 860--914. doi.org/10.1007/s00039-015-0327-6

\bibitem{KanKrotkova} Kan, I. D., Krotkova N. A. (2011). Quantitative generalizations of Niederreiter's results on continued fractions. \textit{Chebyshevsk\u{\i}i Sb.} 12(1): 100-119.

\bibitem{komatsu} Komatsu, T. (2005). On a Zaremba's conjecture for powers. \textit{Sarajevo J. Math.} 1(13): 9--13.

\bibitem{kontorovich} Kontorovich, A. (2013). From Apollonius to Zaremba: local-global phenomena in thin orbits. \textit{Bull. Amer. Math. Soc. (N.S.)} 50(2): 187--228. doi.org/10.1090/S0273-0979-2013-01402-2

\bibitem{mendesfrance} Mendès France, M. (1973). Sur les fractions continues limitées. \textit{Acta Arith.} 23(2): 207--215. doi.org/10.4064/aa-23-2-207-215 

\bibitem{Niederreiter1986} Niederreiter, H. (1986). Dyadic fractions with small partial quotients. \textit{Monatsh. Math.} 101: 309--315. doi.org/10.1007/BF01559394

\bibitem{vanderPoortenShallit} van der Poorten, A. J., Shallit, J. (1992). Folded continued fractions. \textit{J. Number Theory.} 40(2): 237--250. doi.org/10.1016/0022-314x(92)90042-n

\bibitem{YodphotongLaohakosol2002} Yodphotong, M., Laohakosol, V. (2002). Proofs of Zaremba's Conjecture for powers of 6. In: \textit{Proceedings of the International Conference on Algebra and Its Applications (ICAA 2002) (Bangkok)}. pp. 278--282. 

\bibitem{Zaremba} Zaremba, S. K. (1972). La Méthode des “Bons Treillis” pour le Calcul des Intégrales Multiples. In: \textit{Applications of Number Theory to Numerical Analysis}. (Proc. Sympos., Univ. Montreal, Montreal). pp. 39--119.

\end{thebibliography}
\end{document}